%% file: rpb200.tex
\documentclass{article}

\begin{document}
\newcommand{\lbrk}{\linebreak[0]}
\newcommand{\pp}{{\em{pp.}}}
\bibliographystyle{plain}

\title{Factorizations of Cunningham numbers with bases $13$ to $99$:\\
Millennium edition%
\thanks{
\hbox{Incorporating {\em
Factorizations of $a^n \pm 1$, $13 \le a < 100$: Update 3.}
This version corrected 18 March 2001.}
}
\thanks{
\hbox{Copyright \copyright\ 2000, the authors.\hspace{10em}}
\hspace*{\fill} \hbox{rpb200 typeset using \LaTeX}
}}

\author{\\[20pt]
	Richard P.\ Brent\\
	Oxford University Computing Laboratory\\
	Wolfson Building, Parks Road\\
	Oxford OX1 3QD, UK\\[-2pt]
	{\tt\small rpb@comlab.ox.ac.uk}
\and
	\\[20pt]
	Peter L.\ Montgomery\\
	780 Las Colindas Road\\
	San Rafael, CA 94903--2346\\
	USA\\[-2pt]
	{\tt\small pmontgom@cwi.nl}
\and
	\\[20pt]
	Herman J.\ J.\ te Riele\\
	Centrum voor Wiskunde en Informatica\\
	Kruislaan 413\\
	1098 SJ Amsterdam\\
	The Netherlands\\[-2pt]
	{\tt\small herman@cwi.nl}\\[25pt]
	}

\date{with the assistance of\\[15pt]
Henk Boender, Stephania Cavallar, Conrad Curry,\\
Bruce Dodson, Jens Franke, Joseph Leherbauer,\\
George Sassoon and Robert Silverman
\\[20pt]
PRG TR-14-00\\[20pt]
31 December 2000\\}

\maketitle
\vspace{-15pt}
\pagenumbering{roman}			
\thispagestyle{empty}			

\vspace{\fill}
\begin{abstract}

This Report updates the tables of factorizations
of $a^n \pm 1$ for $13 \le a < 100$, previously published
as CWI Report NM-R9212 (June 1992)
and updated in CWI Report NM-R9419 (Update~1, September 1994)
and CWI Report NM-R9609 (Update~2, March 1996).
A total of 951 
new entries in the tables are given here.
The factorizations are now complete for $n < 76$, 
and there are no composite cofactors smaller than
$10^{102}$.\\[10pt] 
This ``Millennium edition'' gives the complete tables incorporating
all updates. A file containing only the new updates, a file containing
factorizations for an extended table range, and a file of factors,
are all available on the internet.\\

{\em 1991
Mathematics Subject Classification.}
Primary 11A25; Secondary 11-04\\
{\em Key words and phrases.}
Cunningham numbers,
elliptic curve method,
factor tables,
number field sieve,
quadratic sieve,
ECM, GNFS, MPQS, NFS, PMPQS, PPMPQS, SNFS
\end{abstract}

\vspace*{\fill}
\pagebreak[4]
\section{Introduction}
\label{sec:Intro}

For many years there has been an interest in the prime factors of numbers
of the form $a^n \pm 1$, where~$a$ is a small integer (the
{\em base}) and~$n$ is a positive exponent.
Such numbers often arise.  For example, if~$a$ is prime then
there is a finite field~$F$ with~$a^n$ elements, and the multiplicative
group of~$F$ has $a^n - 1$ elements.  Also, for prime~$a$
the sum of divisors of~$a^n$
is $\sigma(a^n) = (a^{n+1} - 1)/(a - 1)$.  Numbers of the form
$a^n + 1$ arise as factors of $a^{2n} - 1$ and in other ways.

An extensive table of factors of $a^n \pm 1$ for $a \le 12$ has been published
by Brillhart {\em et~al}~\cite{Bri88a}.  The computation of these tables
is referred to as the {\em Cunningham Project}
in recognition of the pioneering
computations of Cunningham and Woodall~\cite{Cun25}.  For a history, see the
Introduction in~\cite{Bri88a}.

The tables~\cite{Bri88a} are limited to $a \le 12$,
but many applications require larger bases.
In June 1992 tables covering the range $13 \le a < 100$
were published~\cite{rpb134}.
The exponents~$n$ satisfied $a^n < 10^{255}$ if $a < 30$, and
$n \le 100$ if $a \ge 30$. An update~\cite{rpb134u1} containing
780	
new factorizations (with the same limits for $a$ and $n$)
was published in September 1994, and a second update~\cite{rpb134u2}
containing 760 new factorizations was published in March 1996.
These factorizations are now incorporated in the
{\em Magma} package~\cite{Bos96}.

Since the second update~\cite{rpb134u2},
many new factors have been found.
The factorizations are now complete for $n \le 75$, 
and there are no composite cofactors
with fewer than~$103$ digits\footnote{``digits'' always means
``decimal digits''.}. 
This report includes all the new (complete or partial) factorizations
found from the publication of~\cite{rpb134u2} to 31 December 2000.
Altogether, 
951 new (complete or partial) 
factorizations are listed,
involving 1098 new  
factors\footnote{Here
and elsewhere we do not count large factors which
are obtained by division by other factors.}.
Table~\ref{Tab:Progress} summarizes progress since the publication
of the original tables~\cite{rpb134}. ``Update~3'' refers to this Report.

\begin{table}[h]
\caption{Statistics regarding the Tables and Updates}
~\\[-5pt]	
\label{Tab:Progress}
\centerline{
\begin{tabular}{|c|c|c|c|c|c|} \hline
Tables	 & Date		& Smallest 	& Complete to & Total  \\
	 &		& composite	& exponent    & entries\\
									\hline
Original & June 1992	& 81 digits	& 46	      & 13882  \\
Update~1 & Sept.~1994	& 87 digits     & 58	      & 780    \\
Update~2 & March 1996	& 95 digits	& 66	      & 760    \\
Update~3 & Dec.~2000    & 103 digits    & 75	      & 951	\\
									\hline
\end{tabular}
}
\end{table}

Table~\ref{Tab:Distribution} shows the number of prime factors of different
sizes found for Updates~1--3 (excluding large factors obtained
by division).  The median sizes are 26~digits for Update~1,
29--30~digits for Update~2, 
and 33~digits for Update~3.
The {\em smallest} new factor is 20 digits for Update~3
(compare 14 digits for Update~2). We would be surprised
if many factors of less than 25 digits are still to be found.
The {\em largest} new penultimate factor is 60 digits for Update~3
(compare 56 digits for Update~2).

\begin{table}[ht]
\caption{Distribution of Factors}
~\\[-5pt]	
\label{Tab:Distribution}
\centerline{
\begin{tabular}{|c|c|c|c|} \hline
Digits	 & Update 1	& Update 2 	& Update 3	\\ \hline
10--14	& 0		& 1	& 0  \\
15--19	& 17		& 23	& 0  \\
20--24	& 333		& 144	& 24 \\
25--29	& 329		& 273	& 242\\
30--34	& 154		& 197	& 322\\
35--39	& 72		& 99	& 181\\
40--44	& 44		& 89	& 134\\
45--49	& 9		& 39	& 107\\
50--54	& 0		& 14	& 63\\
55--59	& 1		& 3	& 15\\
60--64	& 0		& 0	& 10\\ \hline
Total	& 959		& 882	& 1098 \\	\hline
\end{tabular}
}
\end{table}

\vspace*{\fill}
\pagebreak[4]
\section{Format of the Tables}

The format of the Tables is the same as in~\cite{rpb134}.
For each base~$a$, not a perfect power, in the range $13 \le a < 100$,
we give two separate tables~--

\begin{itemize}
\item[] {\em Table a}$-$: factorizations of $a^n - 1$, $n$ odd.\\[-20pt]
\item[] {\em Table a}$+$: factorizations of $a^n + 1$.
\end{itemize}

\noindent The exponent ranges are as in~\cite{rpb134}~--
\begin{itemize}
\item[] $13 \le a < 30$, exponents $n$
	such that $a^n < 10^{255}$.			\\[-20pt]
\item[] $30 \le a < 100$, exponents $n \le 100$.
\end{itemize}

The entries are similar in format to those of the ``short'' tables
in~\cite{Bri88a}.
All known factors, including algebraic and
Aurifeuillian~\cite{rpb135} factors, are listed.
Factors which are given as decimal numbers are primes.
Exponents are indicated by a hat~({\char94}), for example
``2{\char94}3'' means~$2^3$.
Multiplication is indicated by a period~(.),
for example $3^3 + 1 = 2^27$ is written as ``2{\char94}2.7''.
A period at the end of a line implies that the factorization is
continued on the next line.
An underscore (\_) at the end of a line means that a (large) factor
is continued on the next line (see, for example, the entry for $19^{177}-1$).

The largest factor of $a^n \pm 1$ may be found by division
by the smaller factors.  Thus, such factors are abbreviated.
The notation $p_{xy}$ or ``{\tt{pxy}}''
means a prime factor of {\tt{xy}} digits.
For example, the prime 1238926361552897 might be abbreviated as~{\tt{p16}}.
Similarly, the notation $c_{xy}$ or ``{\tt{cxy}}''
means a composite number of {\tt{xy}} digits.

\section{Availability of Tables, Updates and Factors}
\label{sec:avail}

The complete tables for the original table range are given below.
The changes since Update~2 are available by anonymous ftp from
{\tt ftp://ftp.comlab.ox.ac.uk/{\lbrk}pub/{\lbrk}Documents/%
{\lbrk}techpapers/{\lbrk}Richard.Brent/{\lbrk}rpb134u3.txt.gz}
(a compressed text file). This file includes comments on the person
and method responsible for finding each factor (if there is no attribution,
the factor was found by one of the authors).

The complete tables incorporating Updates~1--3, and a list of factors,
are also available online:
see
{\tt http://www.comlab.ox.ac.uk/{\lbrk}oucl/{\lbrk}work/%
{\lbrk}richard.brent/{\lbrk}factors.html}
The restriction
$n \le 100$ for bases $a \ge 30$ has been relaxed for the online tables;
we only require $a^n < 10^{255}$.
For this extended table range, the smallest composite has 102 digits.

\section{Factorization Methods}

Since Update 2 we have attempted to
factor the remaining composite numbers in the tables by using the
{\em elliptic curve method} (ECM).
Sometimes ECM is successful in finding one or more factors.
If the factorization can not be completed by ECM, but the remaining composite
part is sufficiently small, we use the
{\em multiple polynomial quadratic sieve} (MPQS) method
to complete the factorization. In some cases we prefer to use
the {\em number field sieve} (NFS) if it is predicted to be
faster than MPQS\footnote{The choice depends upon
the size of the known non-algebraic factors of the number $a^n \pm 1$.
We normally use the {\em special} number field sieve (SNFS),
but in at least one case
($17^{186}+1$) the {\em general} number field sieve (GNFS)
was used (this factorization was contribution by
Couveignes, Granboulan, Hoogvorst and Nguyen).}.

We do not describe ECM, MPQS or NFS here.
The reader should refer to~\cite{Len87,Mon87,Mon94}
for a general description of ECM,
to~\cite{Boe95,Sil87} for MPQS, and to~\cite{Len93a,Hui95,Mur99} for NFS.
A recent survey is~\cite{rpb196}.
The particular implementations of ECM by Brent and Montgomery are
described in~\cite{rpb161,Mon92}.

Table~\ref{Tab:counts} shows the number of factors found by
several methods in the preparation of Updates~1--3.
For ECM and MPQS these only include penultimate factors
of at least 30 digits.
An increase in the use of SNFS
and decline in the use of Pollard's $p \pm 1$ methods~\cite{Rie94} is evident.
There is also a marked increase in the number of large (at least 30-digit)
factors found by ECM. Most of the new factors found by MPQS and SNFS are large
because these methods are only used after ECM has been tried.
In fact, since Update~2, MPQS and SNFS did not find any factor with
less than thirty digits, because such factors had already been found by ECM.
The largest factor found by ECM was a 52-digit factor of
$96^{98}+1$ (see~\cite{champ40}).

\begin{table}[ht]
\caption{Factors Found by Different Methods}
~\\[-5pt]	
\label{Tab:counts}
\centerline{
\begin{tabular}{|c|c|c|c|c|c|} \hline
	 & Pollard & Pollard & ECM     & MPQS    & NFS \\
	 & $p-1$   & $p+1$   & $(30$D$+)$ & $(30$D$+)$ & \\ \hline
Update~1 &  38     &   16    & 	69     & 157      & 37   \\
Update~2 &   0     &    3    & 	151    & 155    & 136  \\
Update~3 &   0     &    3    & 	423    & 129    & 279 \\ \hline 
\end{tabular}
}
\end{table}

\vspace*{\fill}
\pagebreak[4]
\section{First Holes}

A ``first hole'' is the first composite number occurring in a table.
Thus, each table of factorizations is complete up to, but not including,
its first hole. Table~\ref{Tab:Holes} lists the
exponents of the current first holes
for $2 \le a < 100$ (the range $2 \le a \le 12$ is included for
the sake of comparison).
For example, the first holes in the tables for $a = 17$ occur for
exponents 137 and 118. In fact, first holes such as
$17^{118}+1 = 2 \cdot 5 \cdot 29 \cdot 7789 \cdot c_{139}$
are good candidates for factorization by SNFS.

\begin{table}[htb]
\caption{Exponents of First Holes for $2 \le a \le 99$}
~\\[-5pt]	
\label{Tab:Holes}
\centerline{
\begin{tabular}{|c|c|c||c|c|c||c|c|c||c|c|c|} \hline
$a$	& $-$& $+$&
$a$	& $-$& $+$&
$a$	& $-$& $+$&
$a$	& $-$& $+$\\ \hline
2 & 641 & 617 &
3 & 379 & 382 &
5 & 307 & 283 &
6 & 251 & 232 \\
7 & 227 & 214 &
10 & 197 & 223 &
11 & 191 & 181 &
12 & 173 & 172 \\
13 & 161 & 151 &
14 & 149 & 134 &
15 & 127 & 122 &
17 & 137 & 118 \\
18 & 131 & 121 &
19 & 155 & 113 &
20 & 149 & 106 &
21 & 125 & 128 \\
22 & 103 & 116 &
23 & 101 & 101 &
24 & 101 & 107 &
26 & 107 & 103 \\
28 & 103 & 106 &
29 & 101 & 112 &
30 & 103 & 103 &
31 & 97 & 113\\
33 & 103& 89 &
34 & 103& 101&
35 & 97 & 103&
37 & 89 & 97 \\
38 & 101& 86 &
39 & 89 & 89 &
40 & 97 & 97 &
41 & 101& 89 \\
42 & 115& 86 &
43 & 101& 89 &
44 & 103& 94 &
45 & 83 & 92 \\
46 & 101& 82 &
47 & 89 & 86 &
48 & 107& 94 &
50 & 89 & 97 \\
51 & 97 & 83 &
52 & 83 & 82 &
53 & 89 & 88 &
54 & 107& 79 \\
55 & 107& 86 &
56 & 83 & 79 &
57 & 83 & 79 &
58 & 83 & 76 \\
59 & 97 & 79 &
60 & 79 & 86 &
61 & 79 & 94 &
62 & 79 & 82 \\
63 & 83 & 83 &
65 & 79 & 79 &
66 & 97 & 86 &
67 & 83 & 82 \\
68 & 79 & 76 &
69 & 83 & 83 &
70 & 83 & 89 &
71 & 89 & 76 \\
72 & 79 & 83 &
73 & 79 & 83 &
74 & 89 & 82 &
75 & 79 & 79 \\
76 & 103& 79 &
77 & 89 & 86 &
78 & 79 & 79 &
79 & 97 & 83 \\
80 & 83 & 82 &
82 & 83 & 79 &
83 & 79 & 82 &
84 & 79 & 76 \\
85 & 83 & 76 &
86 & 79 & 79 &
87 & 97 & 83 &
88 & 79 & 92 \\
89 & 83 & 83 &
90 & 83 & 79 &
91 & 83 & 92 &
92 & 83 & 82 \\
93 & 79 & 82 &
94 & 79 & 76 &
95 & 79 & 79 &
96 & 83 & 79 \\
97 & 83 & 82 &
98 & 95 & 79 &
99 & 89 & 88 &
   &    &    \\ \hline
\end{tabular}
}
\end{table}

\section{Probable Primes}

Numbers listed as prime have not in all cases been rigorously
proved to be prime; they may merely have passed a probabilistic
primality test~\cite{Knu81}.  There is a positive but extremely small
probability that a composite number
will pass such a test and be mistaken for a prime.
In applications where it is essential for primality to be
proven rigorously, one should apply an algorithm such as
Morain's elliptic curve primality test~\cite{Atk93b,Mor90b},
which can easily prove
or disprove the primality of numbers of the size considered here.

\pagebreak[4]
\section*{Acknowledgements}

We thank the people mentioned on the title
page for their assistance since Update~2.
In particular, Robert Silverman helped eliminate many first holes.

Thanks to Franz-Dieter Berger, Anders Bj\"orn, Henk Boender, Arjen Bota,
Lars Brunjes, John Cannon, Stephania Cavallar, Graeme Cohen, Conrad Curry,
Nicolas Daminelli, Thomas Denny, Bruce \hbox{Dodson},
Harvey Dubner, Robert Dubner,
Will Edgington, Jens Franke, Torbj\"orn Granlund,
Tom Hill, 
Uwe \hbox{Hollerbach}, Rob Hooft, Marije Huizing,
Wilfrid Keller, Yuji Kida, Hidenori Kuwakado,
Samuli \hbox{Larvala}, Joseph Leherbauer, Arjen Lenstra,
Paul Leyland, Walter Lioen, Hisanori Mishima, Peter Moore,
\hbox{Mitsuo} Morimoto, Andreas M\"uller, Henrik Olsen, John Pollard,
Hans Riesel, David Rusin, George \hbox{Sassoon} (Mullfac), Robert Silverman,
Thomas Sosnowski, Allan Steel, Andrew Steward, Samuel Wagstaff, Georg Wambach,
Damian Weber, Michael Wiener, Hugh Williams, Aiichi Yamasaki, Paul Zimmermann,
those mentioned in a previous footnote,
and anyone inadvertently omitted from this list, for sending us
factors or assisting in other ways.

The Australian National University Supercomputer Facility
provided computer time to run the first author's ECM programs
on a Fujitsu VP~2200/10 vector processor and an SGI Power Challenge.
The ANU-Fujitsu CAP Research Project provided time on a Fujitsu AP~1000.

The Oxford Supercomputing Centre provided time to run the first author's
ECM programs on {\em Oscar}, an SGI Origin 2000.

The Dutch National Computing Facilities Foundation
provided computer time
to run the programs of Boender, Elkenbracht-Huizing,
Montgomery, and te~Riele.

Many runs with ECM, MPQS and SNFS were carried out on
workstations at CWI and
Leiden University.
We are grateful to the workstation
``owners'' at CWI and Leiden University for letting us use their idle
cycles for this project.

\pagebreak
\vspace*{100pt}

\section*{\hspace*{\fill}\hbox{Part II -- The Tables}\hspace*{\fill}}
\pagenumbering{arabic}			
\vspace*{30pt}

There have been many improvements
since the original tables were published in June 1992.
The original tables included $3093$ composites (i.e.~incomplete
factorizations); this number has been reduced to $1027$.
Following are the complete updated tables
with the same range as the original table.
See \S\ref{sec:avail} for information regarding online versions
of these tables (with an extended table range) and associated factors.
\vspace*{\fill}

\input rpb200t.tex		

\end{document}

%% file: rpb200t.tex
\pagebreak
\pagebreak
\pagebreak
\pagebreak
\pagebreak
\pagebreak
\pagebreak
\pagebreak
\pagebreak
\pagebreak
\pagebreak
\pagebreak
\pagebreak
\pagebreak
\pagebreak
\pagebreak
\pagebreak
\pagebreak
\pagebreak
\pagebreak
\pagebreak
\pagebreak
\pagebreak
\pagebreak
\pagebreak
\pagebreak
\pagebreak
\pagebreak
\pagebreak
\pagebreak
\pagebreak
\pagebreak
\pagebreak
\pagebreak
\pagebreak
\pagebreak
\pagebreak
\pagebreak
\pagebreak
\pagebreak
\pagebreak
\pagebreak
\pagebreak
\pagebreak
\pagebreak
\pagebreak
\pagebreak
\pagebreak
\pagebreak
\pagebreak
\pagebreak
\pagebreak
\pagebreak
\pagebreak
\pagebreak
\pagebreak
\pagebreak
\pagebreak
\pagebreak
\pagebreak
\pagebreak
\pagebreak
\pagebreak
\pagebreak
\pagebreak
\pagebreak
\pagebreak
\pagebreak
\pagebreak
\pagebreak
\pagebreak
\pagebreak
\pagebreak
\pagebreak
\pagebreak
\pagebreak
\pagebreak
\pagebreak
\pagebreak
\pagebreak
\pagebreak
\pagebreak
\pagebreak
\pagebreak
\pagebreak
\pagebreak
\pagebreak
\pagebreak
\pagebreak
\pagebreak
\pagebreak
\pagebreak
\pagebreak
\pagebreak
\pagebreak
\pagebreak
\pagebreak
\pagebreak
\pagebreak
\pagebreak
\pagebreak
\pagebreak
\pagebreak
\pagebreak
\pagebreak
\pagebreak
\pagebreak
\pagebreak
\pagebreak
\pagebreak
\pagebreak
\pagebreak
\pagebreak
\pagebreak
\pagebreak
\pagebreak
\pagebreak
\pagebreak
\pagebreak
\pagebreak
\pagebreak
\pagebreak
\pagebreak
\pagebreak
\pagebreak
\pagebreak
\pagebreak
\pagebreak
\pagebreak
\pagebreak
\pagebreak
\pagebreak
\pagebreak
\pagebreak
\pagebreak
\pagebreak
\pagebreak
\pagebreak
\pagebreak
\pagebreak
\pagebreak
\pagebreak
\pagebreak
\pagebreak
\pagebreak
\pagebreak
\pagebreak
\pagebreak
\pagebreak
\pagebreak
\pagebreak
\pagebreak
\pagebreak
\pagebreak
\pagebreak
\pagebreak
\pagebreak
\pagebreak